\newcommand{\ind}{\mbox{{\rm index} \,}}
\newcommand{\coker}{\mbox{{\rm coker} \,}}
\newcommand{\tr}{{\rm tr\,}}
\newcommand{\range}{{\rm range\,}}
\newtheorem{theorem}{Theorem}

\documentstyle[12pt]{article}
\begin{document}

\title{Note on the Curvature and Index of Almost Unitary Contraction Operator
 }
\author{R.N. Levy \\
Department of Mathematics and Informatics, Sofia University  \\ 5
J.Bourchier Blvd., Sofia 1164, Bulgaria\\
e-mail:levy@fmi.uni-sofia.bg
\date{November 29, 2000}
}
 \maketitle
\begin{abstract}
In the recent preprint \cite{par} S. Parrott proves the equality between
the Arveson`s curvature and the Fredholm index of a "pure"
contraction with finite defect numbers. In the present note one
derives a similar formula in the "non-pure" case.
\end{abstract}

The notions of $d$-{\it contraction} $T = (T_1, T_2, \ldots ,
T_d)$ and its {\it curvature} was introduced by W. Arveson in a
series of papers (see  \cite{arv1}, \cite{arv2},  and
\cite{arv3}). In the case of a single contraction ($d=1$) the
curvature is thoroughly investigated in the paper of Parrott
\cite{par}. Namely, let $T$ be a contraction operator on a Hilbert
space $H$, and suppose that $\Delta_T:=\sqrt{I-TT^*}$ has finite
rank. Parrott shows that the curvature $K(T)$ of $T$ can be
defined on three equivalent ways:

\begin{eqnarray}
\nonumber
 K(T)&=& \int_{|z| = 1}  dz \, \lim_{r \uparrow 1} \,
(1-r^2) \, \tr (\Delta_T (I-rzT^*)^{-1} (I-r \bar{z} T)^{-1}
\Delta_T) \\ \nonumber
 &=&\lim_{n \to \infty} \frac{\tr( I - T^n {T^*}^n)}{n} \\
 \nonumber
 &=& \lim_{n \to \infty} \tr ({T^*}^n T^n (I - TT^*))\quad .
 \end{eqnarray}

In the papers cited above Arveson introduces the notion of "pure"
$d$ - contraction. In the case of a single contraction $T$ this
reduces to the condition that $T$ belongs to the class
$C_{\cdot,0}$ , i.e for any $h \in H$ we have $T^{*n}h \to 0$ as $n
\to \infty$. For a single pure contraction Parrott proves that
$$ K(T)\; = \; - \; \ind T$$

We prove a similar formula in the general "non-pure" case. Note
that in the "pure" case (as it is noted in the Parrott's paper)
the assumption of the finiteness of the rank of the operator $I- T
T^*$ implies that the rank of $I - T^* T$ is also finite (and not
exceeding the rank of $I- T T^*$). In the "non-pure" case we need
to postulate this; more generally, recall that an operator $T$ in
Hilbert space is called {\it almost unitary} if both $I - T^*T$
and $I-TT^*$ are trace-class operators.

\begin{theorem}
Let $T$ be almost unitary contraction. Then
$$\ind T \, = \, K(T^*) \, - \, K(T)$$
\end{theorem}

{\bf Proof.} Denote $a_n(T):= \tr ({T^*}^n T^n (I - TT^*))$. Then
$$a_n(T)= \tr ( T^n (I - TT^*){T^*}^n) = \tr \left( T^n {T^*}^n -
T^{n+1} {T^*}^{n+1}\right)$$ On the other side, $$a_n(T)= \tr
\left({T^*}^n T^n - {T^*}^{n+1} T^{n+1} + \left[ T^*, {T^*}^n
T^{n+1} \right]\right) = a_n(T^*) + b_n(T)$$ where $b_n(T) := \tr
\left[ T^*, {T^*}^n  T^{n+1} \right]$.

Since $K(T) = \lim_{n \to \infty} a_n(T)$ , $K(T^*) = \lim_{n \to
\infty} a_n(T^*)$ , the proof will be completed if one can show
that $b_n(T) = - \,\ind T$ for $n \geq 0$.

For this, note first that since $T$ and $T^*$ commute modulo trace
class operators, the trace of the commutator above does not depend
on the order of factors in the second term, and therefore $b_n(T)
= \tr \left[ T^*, T(T^* T)^n  \right].$ It is easy to see that
$b_0(T)=\left[ T^*,T\right]= -\ind T$.\footnote{ This well-known
fact can be proven in the following way: denote by $S$ the unique
operator determined by the conditions $\ker S = (\range T)^\bot\,
, \, \range S = (\ker T)^\bot$ , and $TSx = x$ for any $x \in
\range T$. Then $\ind T = \dim\ker T - \dim \coker T = \tr (I-
ST)\, - \, \tr(I-TS)= \tr[T,S]$. On the other hand, $S-T^*$ is
trace-class operator (both $S$ and $T^*$ are inverse to $T$ modulo
trace-class operators) and therefore $\tr[T,S]=\tr[T,T^*]$.} On
the other hand, since $I-T^*T$ is trace-class, then for any $n
\geq 0$ one has $$b_n(T) \,-\, b_{n+1}(T) \,=\, \tr
\left[T^*,T(T^*T)^n\left( I-T^*T\right) \right]=0,$$which
completes the proof.

\medskip

{\bf Remark 1.} If $T$ is pure, then $K(T^*)=0$, and one recovers
the result from \cite{par} cited above.

\medskip

{\bf Remark 2.} An alternative way to compute $b_n(T)$ is to use
the Helton - Howe formula for traces of commutators. In our case
this gives no profits, but one may expect that this approach could
be useful in the multidimensional case.

\medskip

The author thanks Stephen Parrott for his useful remarks.

\end{document}